\documentclass[13pt, reqno]{amsart}
\usepackage{amsmath}
\usepackage{amscd}
\usepackage{amsthm}
\usepackage{amssymb} \usepackage{latexsym}
\usepackage{eufrak}
\usepackage{euscript}
\usepackage{epsfig}
\usepackage{graphics}
\usepackage{array}
\usepackage{enumerate}
\usepackage{color}
\usepackage{wasysym}
\usepackage{pdfsync}
\usepackage{stmaryrd}

\newcommand{\bel}[1]{\begin{equation}\label{#1}}

\newcommand{\be}{\begin{equation}}

\newcommand{\dist}{\rho}
\newcommand{\ba}{\begin{eqnarray}}
\newcommand{\ea}{\end{eqnarray}}

\newcommand{\qe}{\end{equation}}
\newcommand{\R}{{\mathbb R}}
\newcommand{\N}{{\mathbb N}}
\newcommand{\Z}{{\mathbb Z}}

\newcommand{\wt}{\widetilde}

\newcommand{\eg}{\begin{example}}
\newcommand{\egd}{\end{example}}
\newcommand{\tm}{\begin{thm}}
\newcommand{\tmd}{\end{thm}}
\newcommand{\co}{\begin{coro}}
\newcommand{\cod}{\end{coro}}
\newcommand{\enu}{\begin{enumerate}}
\newcommand{\enud}{\end{enumerate}}
\newcommand{\rmk}{\begin{rem}}
\newcommand{\rmkd}{\end{rem}}

\theoremstyle{theorem}
\newtheorem{thm}{Theorem}[section]
\newtheorem{prop}{Proposition}[section]
\theoremstyle{example}
\newtheorem{example}{Example}[section]
\theoremstyle{corollary}
\newtheorem{coro}{Corollary}[section]
\theoremstyle{lemma}

\theoremstyle{definition}

\theoremstyle{proof}

\theoremstyle{remark}
\newtheorem{rem}{Remark}[section]
\theoremstyle{remark}

\newtheorem{assump}{Assumption}[section]

\begin{document}

\title[Ancient caloric functions on graphs]{Ancient caloric functions on graphs with unbounded Laplacians}

\author{Bobo Hua}
\address{School of Mathematical Sciences, LMNS,
Fudan University, Shanghai 200433, China; Shanghai Center for
Mathematical Sciences, Fudan University, Shanghai 200433,
China.}
\email{bobohua@fudan.edu.cn}

\begin{abstract}
We study ancient solutions of polynomial growth to both continuous-time and discrete-time heat equations on graphs with unbounded Laplacians. We generalize Colding and Minicozzi's theorem \cite{ColdingM19} on manifolds, and the result \cite{Huaancient19} on graphs with normalized Laplacians to the setting of graphs with unbounded Laplacians: For a graph admitting an intrinsic metric, which has polynomial volume growth, the dimension of the space of ancient solutions of polynomial growth is bounded by the dimension of harmonic functions with the same growth up to some factor. \end{abstract}

\maketitle




\section{Introduction}

Let $M$ be a complete, noncompact Riemannian manifold without boundary. For any $k>0,$ we denote by $\mathcal{H}_k(M)$ the space of harmonic functions of polynomial growth with the growth rate at most $k,$ i.e. $f\in \mathcal{H}_k(M)$ if $\Delta f=0$ and there exist $p\in M$ and  a constant $C_f,$ depending on $f,$ such that
$$\sup_{x\in B_R(p)}|f(x)|\leq C_f(1+R)^k,\quad \forall\  R>0,$$ where $B_R(p)$ denotes the ball of radius $R$ centered at $p.$ 

For a Riemannian manifold $M$ with nonnegative Ricci curvature, Yau \cite{YauCPAM75} proved the Liouville theorem that any positive harmonic function on $M$ is constant. Yau conjectured that for any $k>0$ the space $\mathcal{H}_k(M)$ is a finite-dimensional linear space, see e.g. \cite{Yau87EM,Yau93book}. This conjecture was settled in \cite{ColdingMAnnals97}, see also \cite{ColdingMJDG97,ColdingMInv98,ColdingMCPAM98, LiMRL98, CCM95, LiTam89} for related results. 

A natural generalization is to consider ancient solutions, defined on the time interval $(-\infty,0],$ of polynomial growth to heat equations. For a Riemannian manifold $M$ and $k>0,$ we denote by $\mathcal{P}_k(M)$ the space of ancient solutions $u(x,t)$ satisfying that
there exist $p\in M$ and a constant $C_u>0$ such that
$$\sup_{B_R(p)\times [-R^2, 0]}|u|\leq C_u(1+R)^k,\quad \forall\  R>0.$$
Calle \cite{CalleMZ06,Callethesis} initiated the study of dimensional bounds for $\mathcal{P}_k(M).$ For an $n$-dimensional Riemannian manifold $M$ with nonnegative Ricci curvature, Lin and Zhang \cite{LinZhang17} proved that
$$\dim\mathcal{P}_k(M)\leq C(n) k^{n+1}, \quad k\geq 1.$$ Recently, Colding and Minicozzi \cite{ColdingM19} proved the following general result, which yields the improvement of Lin and Zhang's result, $$\dim\mathcal{P}_k(M)\leq C(n) k^{n}, \quad k\geq 1.$$ 
\tm[\cite{ColdingM19}]\label{thm:CM1} If $M$ has polynomial volume growth, i.e. there exist $p\in M$ and constants $C, d_V$ such that $$\mathrm{Vol}(B_R(p))\leq C(1+R)^{d_V},\quad \forall R>0,$$ where $\mathrm{Vol}$ denotes the Riemannian volume, then
$$\dim \mathcal{P}_{2k}(M)\leq (k+1)\dim \mathcal{H}_{2k}(M),\quad \forall k\geq 1.$$
\tmd

Harmonic functions of polynomial growth on graphs have been extensively studied by many authors, e.g. \cite{DelmottePolynomial98,KleinerJAMS10,ShalomTao10,Taoweb2,HJLAGAG13,HornLinLiuYau14,BDKY15,HuaJostLiu15,HuaJostMathZ15,HuaJostTAMS15,MPTY17}.  For ancient solutions of heat equations on graphs, the author \cite{Huaancient19} generalized Colding and Minicozzi's theorem, Theorem~\ref{thm:CM1}, to graphs with normalized Laplacians, see the definition below. In this paper, we extend the result to the more general setting of graphs with (possibly) unbounded Laplacians.

We recall the setting of weighted graphs. Let $(V,E)$ be a locally finite, simple, undirected graph. Two vertices $x,y$ are called neighbours, denoted by $x\sim y$, if there is an edge connecting $x$ and $y,$ i.e. $\{x,y\}\in E.$ A graph is called connected if for all $x,y\in V,$ there are vertices $z_i,$ $0\leq i\leq n$, such that $x = z_0 \sim. . . \sim z_n = y.$ 
We always assume that the graph $(V,E)$ is connected. 
Let $$w: E\to (0,\infty),\ \{x,y\}\mapsto w_{xy}=w_{yx},$$ be an edge weight function, and $$m:V\to (0,\infty),\ x\mapsto m_x$$ be a vertex weight function. We denote by $\ell^2(V,m)$ the space of $\ell^2$-summable functions on $V$ with respect to the discrete measure $m.$ 
For any $\Omega\subset V,$ we denote by $$m(\Omega):=\sum_{x\in \Omega}m(\Omega)$$ the $m$-measure of $\Omega.$  We call the quadruple $G=(V,E,m,w)$ a weighted graph. There are no relation between the weights, $w$ and $m,$ a priori.

For a weighted graph $G=(V,E,m,w),$ the Laplace operator $\Delta$ is defined as, for any function $f:V\to \R,$
$$\Delta f (x):= \sum_{y\in V:y\sim x}\frac{w_{xy}}{m_x}\left(f(y)-f(x)\right), \quad\forall x\in V.$$ Note that the Laplacian $\Delta$ depends on the choice of weights $w$ and $m.$ One can show, see e.g. \cite{KellerLenz12}, that $\Delta$ is a bounded operator on $\ell^2(V,m)$ if and only if $$\sup_{x\in V}\frac{\sum_{y\in V:y\sim x}w_{xy}}{m_x}<\infty.$$
Given the edge weight $w,$ if we choose $m_x=\sum_{y\in V:y\sim x}w_{xy}$ for all $x\in V,$ then the corresponding Laplacian is called the \emph{normalized Laplaican}, which is the generator for the simple random walk on $G,$ see e.g. \cite{Woess00}. In this paper, we consider general vertex weights $m,$ for which the Laplacians are possibly unbounded.  

For the analysis on graphs with unbounded Laplacians, Frank, Lenz and Wingert \cite{FLW12} introduced the so-called intrinsic metrics, see e.g. \cite{GHM,KellerLenz12,HKMW13,HKW13JLMS,BHKadv,Huang14pa,Folz14,HShiozawa,HuaKeller14,Bauer-Keller-Wojciechowski15,HuaLin17,BauerHY17,GLLY18} for recent developments. A \emph{(pseudo)metric} is a map $\rho:V\times V\to [0,\infty),$ which is symmetric, satisfies the triangle inequality and $\rho(x,x)=0$ for all $x\in V.$  We denote by
$$s:=\sup_{x\sim y}\dist(x,y)$$ the \emph{jump size} of the metric $\dist.$ For any $R>0,$ we write $B_R(x):=\{y\in V: \dist(y,x)\leq R\}$ for the ball of radius $R$ centered at $x$ with respect to the metric $\rho.$

A metric $\dist$ is called an \emph{intrinsic metric} on $G$ if for any $x\in V,$
\begin{equation}\label{eq:de1}\sum_{y\in V: y\sim x}w_{xy}\dist^2(x,y)\leq m_x.\end{equation} In this paper, we only consider intrinsic metrics satisfying the following assumption.
\begin{assump}\label{ass:1} $\rho$ is an intrinsic metric such that
\begin{enumerate}[(i)]
\item for any $ x\in V, R>0,$ $B_R(x)$ is a finite set, and
\item $\rho$ has finite jump size, i.e. $s<\infty.$
\end{enumerate}
\end{assump}

 A function $f$ on $V$ is called harmonic if $\Delta f=0.$  We denote by $\mathcal{H}_k(G)$ the space of harmonic functions of polynomial growth on $G$ with the growth rate at most $k,$ i.e. $f\in \mathcal{H}_k(G)$ if $f$ is a harmonic function on $V$ and there exist $x_0\in V$ and a constant $C_f$ such that
 $$\sup_{x\in B_R(x_0)}|f(x)|\leq C_f(1+R)^k,\quad \forall R>0.$$
 
We say that $G$ has \emph{polynomial volume growth} with respect to $\rho$ if there are $x_0\in V$ and constants $\alpha,C$ such that \begin{equation}\label{eq:polynomial1}m(B_R(x_0))\leq C(1+R)^{\alpha},\quad \forall R>0.\end{equation}


 In the first part of the paper, we consider ancient solutions of polynomial growth for continuous-time heat equations on graphs.
 Let $\R_-:=(-\infty,0].$ A function $u(x,t)$ on $V\times \R_-$ is called an ancient solution to the (continuous-time) heat equation if \begin{equation}\label{eq:cdef1}\frac{\partial }{\partial t} u(x,t)=\Delta u(x,t),\quad \forall x\in V, t\in \R_-.\end{equation} We denote by $\mathcal{P}_k(G)$  the space of ancient solutions of polynomial growth to the heat equation with the growth rate at most $k,$ i.e. $u\in \mathcal{P}_k(G)$ if $u$ is an ancient solution to the heat equation and there are $x_0\in V$ and a constant $C_u$ such that
 $$\sup_{(x,t)\in B_R(x_0)\times[-R^2, 0]}|u(x,t)|\leq C_u(1+R)^k,\quad \forall R>0.$$ 

The following is the main result of the paper.
\tm\label{thm:main1} Let $G$ be a weighted graph admitting an intrinsic metric satisfying Assumption~\ref{ass:1}. If $G$ has polynomial volume growth, then for all $k\geq 1,$
$$\dim \mathcal{P}_{2k}(G)\leq (k+1)\dim \mathcal{H}_{2k}(G).$$
\tmd
A similar result was obtained for graphs with normalized Laplacians in \cite{Huaancient19}. In this paper, we refined the arguments therein and proved the result for any weighted graph with a (possibly) unbounded Laplacian, which admits an intrinsic metric. In particular, we introduce a modified quantity, defined in \eqref{eq:defht}, to circumvent the difficulties in \cite{Huaancient19}.

In the second part of the paper, we consider ancient solutions of polynomial growth for discrete-time heat equations on graphs.
Let $\Z_-:=\Z\cap (-\infty,0].$ 
A function $v(x,t)$ on $V\times \Z_-$ is called an ancient solution to the discrete-time heat equation if $$v(x,t)-v(x,t-1)=\Delta v(x,t),\quad \forall x\in V, t\in \Z_-.$$
We denote by $\wt{\mathcal{P}}_k(G)$ the space of ancient solutions of polynomial growth to the discrete-time heat equation with the growth rate at most $k,$ i.e. $v\in \wt{\mathcal{P}}_k(G)$ if $v$ is an ancient solution to the discrete-time heat equation and there are $x_0\in V$ and a constant $C_v$ such that
 $$\sup_{(x,t)\in B_R(x_0)\times([-R^2, 0]\cap \Z)}|v(x,t)|\leq C_v(1+R)^k,\quad \forall R>0.$$

\tm\label{thm:main2} Let $G$ be a weighted graph admitting an intrinsic metric satisfying Assumption~\ref{ass:1}. If $G$ has polynomial volume growth, then for all $k\geq 1,$
$$\dim \wt{\mathcal{P}}_{2k}(G)\leq (k+1)\dim \mathcal{H}_{2k}(G).$$
\tmd
Due to the discrete nature of the time in the above theorem, there are some new phenomena for the structure of ancient solutions of polynomial growth, see e.g. Corollary~\ref{coro:cdis1}, compared with Corollary~\ref{coro:c1}.

The paper is organized as follows: In the next section, we recall some basic properties of graphs. In Section~\ref{sec:con}, we prove the parabolic Caccioppoli inequality for the heat equation on graphs, and prove Theorem~\ref{thm:main1}. In Section~\ref{sec:dis}, we study discrete-time heat equations and prove Theorem~\ref{thm:main2}.

In this paper, for simplicity the constants $C$ may change from line to line.

\section{Preliminaries}

Let $G=(V,E,m,w)$ be a weighted graph. 
For convenience, we extend the edge weight function $w$ to $V\times V$ by
setting $w_{xy}=0$ for any pair $(x,y)$ with $x\not\sim y.$ In this way, for a function $f$ on $V$ we may write
$$\sum_{y\in V}w_{xy}f(y)=\sum_{y\in V: y\sim x} w_{xy}f(y).$$ For any $\Omega\subset V,$ we write, for simplicity,
$$\sum_\Omega f:=\sum_{x\in \Omega}f(x)m_x,\ \sum f:=\sum_{x\in V}f(x)m_x,$$ whenever they make sense.
The difference operator $\nabla$ is defined as
$$\nabla_{xy}f=f(y)-f(x),\quad \forall x,y\in V.$$ The following proposition is elementary.
\begin{prop} \begin{equation}\label{eq:basic}\nabla_{xy}(fg)=f(x)\nabla_{xy}g+g(y)\nabla_{xy}f.
\end{equation}
\end{prop}

The ``carr\'e du champ" operator $\Gamma$ is defined as 
$$\Gamma(f)(x)=\frac12\sum_{y\in V}\frac{w_{xy}}{m_x}(f(y)-f(x))^2,\quad x\in V.$$ So that $\Gamma(f)$ is a function on $V,$ which is a discrete analog of $|\nabla f|^2$ for a $C^1$ function $f$ on a manifold.


The following Green's formula is well-known, see e.g. \cite[Theorem~2.1]{Grigoryanbook}. We denote by $C_0(V)$ the set of functions on $V$ of finite support.
\tm For any $f,g:V\to \R,$ if $g\in C_0(V),$ then
\begin{equation}\label{eq:Green}\frac12\sum_{x,y\in V}w_{xy}\nabla_{xy}f\nabla_{xy}g=-\sum_{x\in V}\Delta f(x) g(x)m_x.\end{equation}
\tmd

From now on, we fix $x_0\in V$ as a base vertex. Let $\rho$ be an intrinsic metric satisfying Assumption~\ref{ass:1}. We write $B_R:=B_R(x_0),$ $R>0,$ for simplicity. 
For any $R>0,$ we denote by 
\begin{equation}\label{eq:cut1}\eta_R(x):=\max\left\{0,\min\left\{2-\frac{1}{R}\dist(x,x_0),1\right\}\right\}\end{equation}
 the cut-off function on $B_{2R}.$ One easily sees that $\eta_R$ is supported in $B_{2R},$ and takes the value $1$ on $B_R.$ Moreover, one can show that, for any $x,y\in V$ \begin{equation}\label{eq:lip1}|\nabla_{xy}\eta_R|\leq \frac{1}{R}\dist(x,y).\end{equation} That is, $\eta_R$ is a Lipschitz function with Lipschitz constant at most $\frac1R.$

We consider continuous-time heat equations on graphs, see \eqref{eq:cdef1} for the definition. We denote by $$Q_R:=B_R\times [-R^2,0]$$ the parabolic cylinder of size $R$ at $(x_0, 0).$
For a space-time function $u(x,t)$ on $V\times \R_-,$ 
we denote 
$$\int_{Q_R} u:=\int_{-R^2}^{0}\sum_{x\in B_R}u(x,t)m_xdt.$$ For any $t_0\in \R_-,$
we write
$$\left.\sum_\Omega u\right|_{t=t_0}:=\sum_{x\in \Omega}u(x,t_0)m_x.$$
For a $C^1$ function in time, $u(x,t),$ we write $u_t$ (or $\partial_t u$) for $\frac{\partial }{\partial t} u.$

For discrete-time heat equations on graphs, we write
$$\wt{Q}_R=B_R\times ([-R^2,0]\cap \Z),\quad R>0.$$ 
For a space-time function $u(x,t)$ on $V\times \Z_-,$ 
we denote 
$$\sum_{\wt{Q}_R} u:=\sum_{t=-R^2}^{0}\sum_{x\in B_R}u(x,t)m_x.$$ 

For any function $g:\Z_-\to \R,$ we define the difference operator as
$$D_t g(t_0)=g(t_0)-g(t_0-1),\quad \forall t_0\in \Z_-.$$ 
The function $u:V\times \Z_-\to \R$ is an ancient solution to the discrete-time heat equation  if and only if $$D_t u(x,t)=\Delta u(x,t),\quad \forall (x,t)\in V\times \Z_-.$$

The following propositions are elementary. We omit the proofs here.
\begin{prop}\label{prop:ineq1} For any function $g:\Z_-\to \R$ and any $t\in \Z_-,$
$$D_t(g^2)(t)=2g(t)D_tg(t)-(D_t g(t))^2\leq 2g(t)D_tg(t).$$
\end{prop}

\begin{prop}\label{prop:elem1} For any function $g:\Z_-\to \R,$ any $a,b\in \Z_-,$ $a<b,$
$$\sum_{t=a}^b D_t g= g(b)-g(a-1).$$
\end{prop}

\begin{prop}\label{prop:elem2} For $\{a_i\}_{i=1}^N\subset \R,$ there exists some $j,$ $1\leq j\leq N,$ such that
$$a_j\leq \frac1N\sum_{i=1}^N a_i.$$
\end{prop}


\section{Ancient solutions to continuous-time heat equations}\label{sec:con}
In this section, we study ancient solutions to the heat equation on graphs.
The following is the Caccioppoli type inequality to the heat equation on graphs, see e.g. \cite[(3.12)]{LinZhang17} and \cite[(1.2)]{ColdingM19} for Riemannian manifolds and \cite{Huaancient19} for graphs with normalized Laplacians. 
\tm\label{thm:paracacc} There is a universal constant $C$ such that for any ancient solution $u_t=\Delta u$ and $R\geq s,$
\begin{equation}\label{eq:cacci1} R^2\int_{Q_R}\Gamma(u)+R^4\int_{Q_R}u_t^2\leq C\int_{Q_{9R}} u^2.
\end{equation}
\tmd

\begin{proof} We follow the proof strategy by \cite{ColdingM19}, see also \cite{Huaancient19}. For any $R>0,$ let $\eta=\eta_R,$ where $\eta_R$ is the cut-off function defined in \eqref{eq:cut1}.

We first estimate $\int_{Q_R}\Gamma(u).$ 
Since $u_t=\Delta u,$
\begin{equation}\label{eq:energy1}
\frac{\partial }{\partial t} \left(\sum\eta^2 u^2\right)=\sum 2\eta^2 u\partial_tu=2\sum\eta^2u\Delta u.
\end{equation}

By Green's formula \eqref{eq:Green} and \eqref{eq:basic}, 
\begin{eqnarray*}
&&2\sum\eta^2u\Delta u=
-\sum_{x,y}w_{xy}\nabla_{xy}u\nabla_{xy}(\eta^2 u)\\
&=&-\sum_{x,y}w_{xy}\nabla_{xy}u(\eta^2(x)\nabla_{xy}u+u(y,t)\nabla_{xy}(\eta^2))\\
&=&-\sum_{x,y}w_{xy}|\nabla_{xy}u|^2\eta^2(x)-\sum_{x,y}w_{xy}u(y,t)\nabla_{xy}u\nabla_{xy}\eta(2\eta(x)+\nabla_{xy}\eta)\\
&=&-\sum_{x,y}w_{xy}|\nabla_{xy}u|^2\eta^2(x)-2\sum_{x,y}w_{xy}\eta(x)u(y,t)\nabla_{xy}u\nabla_{xy}\eta-\sum_{x,y}w_{xy}u(y,t)\nabla_{xy}u|\nabla_{xy}\eta|^2.
\end{eqnarray*} 
For the last term on the right hand side of the above inequality, by swapping $x$ and $y,$ the symmetry yields that
\begin{eqnarray*}-\sum_{x,y}w_{xy}u(y,t)\nabla_{xy}u|\nabla_{xy}\eta|^2&=&-\frac{1}{2}\sum_{x,y}w_{xy}(u(y,t)-u(x,t))\nabla_{xy}u|\nabla_{xy}\eta|^2\\
&=&-\frac12\sum_{x,y}w_{xy}|\nabla_{xy}u|^2|\nabla_{xy}\eta|^2\leq 0.
\end{eqnarray*} 
Dropping this term, we get for any $R\geq s,$
\begin{eqnarray}\label{eq:bas11}
\partial_t\left(\sum\eta^2 u^2\right)&\leq&-\sum_{x,y}w_{xy}|\nabla_{xy}u|^2\eta^2(x)-2\sum_{x,y}w_{xy}\eta(x)u(y,t)\nabla_{xy}u\nabla_{xy}\eta\nonumber\\
&\leq&-\sum_{x,y}w_{xy}|\nabla_{xy}u|^2\eta^2(x)+\frac12\sum_{x,y}w_{xy}|\nabla_{xy}u|^2\eta^2(x)+2\sum_{x,y}w_{xy}u^2(y,t)|\nabla_{xy}\eta|^2\nonumber\\
&=&-\frac12\sum_{x,y}w_{xy}|\nabla_{xy}u|^2\eta^2(x)+2\sum_{x,y}w_{xy}u^2(y,t)|\nabla_{xy}\eta|^2\nonumber\\
&\leq&-\sum_{x}\Gamma(u)(x)\eta^2(x)m_x+\frac{2}{R^2}\sum_{x,y\in B_{2R+s}}w_{xy}u^2(y,t)\dist^2(x,y)\nonumber\\
&\leq&-\sum_{}\Gamma(u)\eta^2+\frac{2}{R^2}\sum_{y\in B_{3R}}u^2(y,t)m_y.
\end{eqnarray} where we have used the facts that the jump size of $\dist$ is $s,$ \eqref{eq:lip1} and \eqref{eq:de1}.

Fix $R\geq s.$ For $T>0,$ by integrating the above inequality in time from $-T$ to $0,$ we obtain
\begin{eqnarray}
\int_{-T}^0\sum_{B_R}\Gamma(u)&\leq&\int_{-T}^0\sum\Gamma(u)\eta^2\nonumber\\
&\leq& \frac{2}{R^2}\int_{-T}^0\sum_{B_{3R}}u^2+\left.\sum\eta^2 u^2\right|_{t=-T}\nonumber\\
&\leq &\frac{2}{R^2}\int_{-T}^0\sum_{B_{3R}}u^2+\left.\sum_{B_{2R}} u^2\right|_{t=-T}.\label{eq:est1}
\end{eqnarray}
By the mean value property, there is $T_1\in [R^2,4R^2]$ such that
$$\left.\sum_{B_{2R}}u^2\right|_{t=-T_1}=\frac{1}{3R^2}\int_{-4R^2}^{-R^2}\sum_{x\in B_{2R}}m_xu^2(x,t)dt.$$ By using \eqref{eq:est1} for $T=T_1$ and the above equation, we get
\begin{eqnarray}
\int_{Q_R}\Gamma(u)&\leq& \int_{-T_1}^0\sum_{B_R}\Gamma(u)
\leq \frac{2}{R^2}\int_{-T_1}^0\sum_{B_{3R}}u^2+\left.\sum_{B_{2R}} u^2\right|_{t=-T_1}\nonumber\\
&\leq &\frac{C}{R^2}\int_{Q_{3R}}u^2+\frac{1}{3R^2}\int_{-4R^2}^{-R^2}\sum_{x\in B_{2R}}m_xu^2(x,t)dt\nonumber\\
&\leq&\frac{C}{R^2}\int_{Q_{3R}}u^2.\label{eq:est2}
\end{eqnarray} 

Next we estimate $\int_{Q_R}u_t^2.$ Set \begin{equation}\label{eq:defht}h(t):=\frac12\sum_{x,y}w_{xy}|\nabla_{xy}u(\cdot,t)|^2\eta(x)\eta(y).\end{equation} By differentiating $h(t)$ in time, we get
\begin{equation}\label{eq:energy2}\frac{d}{dt}h(t)=\sum_{x,y}w_{xy}\nabla_{xy}u\nabla_{xy}u_t\eta(x)\eta(y).\end{equation}
By Green's formula \eqref{eq:Green}, we get
 \begin{eqnarray*}
 &&\sum_{x,y}w_{xy}\nabla_{xy}u\nabla_{xy}u_t\eta(x)\eta(y)\\
 &=&\sum_{x,y}w_{xy}\nabla_{xy}u\Big[\nabla_{xy}(u_t\eta^2)-\nabla_{xy}\eta(u_t(x,t)\eta(x)+u_t(y,t)\eta(y))\Big]\\
 &=&-2\sum(\Delta u)u_t\eta^2-2\sum_{x,y}w_{xy}u_t(x,t)\eta(x)\nabla_{xy}u\nabla_{xy}\eta\\
 &=&-2\sum u_t^2\eta^2-2\sum_{x,y}w_{xy}u_t(x,t)\eta(x)\nabla_{xy}u\nabla_{xy}\eta\\
 &=:&I+II,
 \end{eqnarray*} where we have used the symmetrization in the third line.
For the second term $II,$ by \eqref{eq:lip1} and \eqref{eq:de1}, for $R\geq s,$
\begin{eqnarray*}
|II|
&\leq&2\sum_{x,y\in B_{2R+s}}w_{xy}|u_t|(x,t)\eta(x)|\nabla_{xy}u|\frac{\dist(x,y)}{R}\\
&\leq&\frac{1}{R^2}\sum_{x,y\in B_{3R}}w_{xy}|\nabla_{xy} u|^2+\sum_{x,y\in B_{3R}}w_{xy} u_t^2(x,t)\eta^2(x)\dist^2(x,y)\\
&\leq&\frac{2}{R^2}\sum_{B_{3R}}\Gamma(u)+\sum_{}u_t^2\eta^2.
\end{eqnarray*}

Hence 
\begin{eqnarray}\label{eq:bas12}\frac{d}{dt}h(t)\leq-\sum u_t^2\eta^2+\frac{2}{R^2}\sum_{B_{3R}}\Gamma(u).
\end{eqnarray}

Fix $R\geq s.$ For $T>0,$ by integrating the above inequality in time from $-T$ to $0,$ and using the properties of $\eta,$ we have

\begin{eqnarray}
\int_{-T}^0\sum_{B_R}u_t^2&\leq&\int_{-T}^0\sum u_t^2\eta^2\leq \frac{2}{R^2}\int_{-T}^0 \sum_{B_{3R}}\Gamma(u)+
h(-T)\nonumber\\
&\leq &\frac{2}{R^2}\int_{-T}^0 \sum_{B_{3R}}\Gamma(u)+\left.\frac12\sum_{x,y\in B_{2R}}w_{xy}|\nabla_{xy}u|^2\right|_{t=-T}\nonumber\\
&\leq&\frac{2}{R^2}\int_{-T}^0 \sum_{B_{3R}}\Gamma(u)+\left.\sum_{B_{2R}}\Gamma(u)\right|_{t=-T}.
\label{eq:et1}
\end{eqnarray}

By the mean value property, there exists $T_2\in [R^2,4R^2]$ such that
\begin{equation}\label{eq:et3}\left.\sum_{B_{2R}}\Gamma(u)\right|_{t=-T_2}=\frac{1}{3R^2}\int_{-4R^2}^{-R^2}\sum_{x\in B_{2R}}m_x\Gamma(u)(x,t)dt\leq \frac{1}{3R^2}\int_{Q_{2R}}\Gamma(u).\end{equation}
By applying \eqref{eq:et1} for $T=T_2$ and using the above equation, we get
\begin{eqnarray}
\int_{Q_R}u_t^2&\leq&\int_{-T_2}^0\sum_{B_R}u_t^2\leq\frac{2}{R^2}\int_{-T_2}^0 \sum_{B_{3R}}\Gamma(u)+\left.\sum_{B_{2R}}\Gamma(u)\right|_{t=-T_2}\nonumber\\
&\leq &\frac{C}{R^2}\int_{Q_{3R}}\Gamma(u)+\frac{1}{3R^2}\int_{Q_{2R}}\Gamma(u)\nonumber\\
&\leq&\frac{C}{R^2}\int_{Q_{3R}}\Gamma(u)\leq \frac{C}{R^4}\int_{Q_{9R}}u^2,\label{eq:et4}
\end{eqnarray} where we have used \eqref{eq:est2}.

The theorem follows from \eqref{eq:est2} and \eqref{eq:et4}.
\end{proof}

This yields the following corollary.
\co\label{coro:c1} Let $G$ be a weighted graph admitting an intrinsic metric satisfying Assumption~\ref{ass:1}. Suppose that $G$ has polynomial volume growth, i.e. \eqref{eq:polynomial1}, and $u\in \mathcal{P}_{k}(G)$ for some $k>0.$  Then for any $q\in \N,$ $4q>2k+\alpha+2,$
$$\partial_t^q u\equiv 0.$$ In particular, there exist some functions $p_i(x),$ $1\leq i\leq q-1,$ such that
$$u(t,x)=\sum_{i=1}^{q-1}p_i(x)t^{i}.$$
\cod
\begin{proof}
We follow the argument in \cite{ColdingM19}. 
For the first assertion,
since $\Delta$ commutes with $\partial_t,$ for any $i\in\N,$ $\partial_t^i u$ is also an ancient solution of the heat equation.  For any $R\geq s,$ applying Theorem~\ref{thm:paracacc} for $\partial_t^i u$ with $0\leq i\leq q-1,$
we get
\begin{eqnarray*}
\int_{Q_R}|\partial_t^qu|^2&\leq& \frac{C}{R^4}\int_{Q_{9R}} |\partial_t^{q-1}u|^2
\leq \cdots\leq \frac{C(q)}{R^{4q}}\int_{Q_{9^qR}} u^2\\
&\leq& C R^{-4q+2k+\alpha+2},
\end{eqnarray*} where we have used $u\in \mathcal{P}_{k}(G)$ and \eqref{eq:polynomial1}. Therefore, for $4q>2k+\alpha+2,$ by passing to the limit $R\to +\infty,$ we prove that $$\partial_t^q u\equiv 0.$$ This proves the first assertion.

The second assertion follows from the first one.
\end{proof}

Now we prove the main theorem, Theorem~\ref{thm:main1}. 
\begin{proof}[Proof of Theorem~\ref{thm:main1}]
The proof follows verbatim from \cite{ColdingM19}. 
Choose $q\in \N$ such that $4q>4k+\alpha+2.$ By Corollary~\ref{coro:c1}, we have
$$u(x,t)=p_0(x)+p_1(x)t+\cdots+p_{q-1}(x)t^{q-1}.$$ 

Note that $u\in \mathcal{P}_{2k}(G).$ For any fixed $x\in V,$ considering sufficiently negative $t$ in the above equality, we obtain that $$p_i(x)=0, \quad \forall i>k.$$ This yields that
\begin{equation}\label{eq:eqpq1}u(x,t)=p_0(x)+p_1(x)t+\cdots+p_{l}(x)t^{l},
\end{equation} where $l:=\lfloor k\rfloor,$ the greatest integer less than or equal to $k.$ 

We claim that the function $p_i(x),$ $0\leq i\leq l,$ grows polynomially with the growth rate less than or equal to $2k.$ Fix distinct values $-1<t_1<t_2<\cdots<t_l<t_{l+1}=0.$ Set column vectors
$$\beta_j:=(1,t_j,t_j^2,\cdots, t_j^l)^T, \quad 1\leq j\leq l+1.$$ We define the matrix
$$B:=(\beta_1,\beta_2,\cdots,\beta_{l+1}).$$ Note that $\det B$ is a Vandermonde determinant, which yields that $\{\beta_j\}_{j=1}^{l+1}$ are linear independent in $\R^{l+1}.$ Let $\{e_i\}_{i=1}^{l+1}$ be the standard basis of $\R^{l+1}.$ We write $b_{j}^i$ for the $ji$-th entry of $B^{-1},$ which implies that
$$e_i=\sum_{j=1}^{l+1} b_j^i\beta_j.$$ Using this fact and \eqref{eq:eqpq1}, we get
$$p_i(x)=\sum_{j=1}^{l+1}b_j^i u(x,t_j).$$ Since $u(x,t_j),$ $1\leq j\leq l+1,$ grows polynomially with the growth rate less than or equal to $2k,$ so does $p_i.$ This proves the claim.

Since $\partial_t u=\Delta u,$ by \eqref{eq:eqpq1}, 
\begin{equation*} \Delta p_l=0,\ \Delta p_i=(i+1)p_{i+1}, \quad 0\leq i\leq l-1.
\end{equation*} Hence we get a linear map
\begin{eqnarray*}\Psi_0:&&\mathcal{P}_{2k}(G)\to \mathcal{H}_{2k}(G)\\
&&\ \  \quad\quad u\mapsto p_l.
\end{eqnarray*} Let $\mathcal{K}_0:=\mathrm{Ker}(\Psi_0).$ It follows that
$$\dim \mathcal{P}_{2k}(G)\leq \dim \mathcal{K}_0+\dim \mathcal{H}_{2k}(G).$$
To estimate $\dim \mathcal{K}_0,$ we note that
for any $u\in \mathcal{K}_0,$ $$p_l=0, \ \Delta p_{l-1}=0.$$ Hence we have a linear map
\begin{eqnarray*}\Psi_1:&&\mathcal{K}_0\to \mathcal{H}_{2k}(G)\\
&&\ \  u\mapsto p_{l-1}.
\end{eqnarray*} Let $\mathcal{K}_1:=\mathrm{Ker}(\Psi_1).$  This yields that
$$\dim \mathcal{K}_0\leq \dim \mathcal{K}_1+\dim \mathcal{H}_{2k}(G).$$ Repeating this $l+1$ times, we prove that
\begin{equation*}\dim \mathcal{P}_{2k}(G)\leq (l+1) \mathcal{H}_{2k}(G)\leq (k+1) \mathcal{H}_{2k}(G).
\end{equation*} This proves the theorem.

\end{proof}

\section{Ancient solutions to discrete-time heat equations}\label{sec:dis}
In this section, we study ancient solutions to the discrete-time heat equation on graphs.

The following is the Caccioppoli type inequality to the discrete-time heat equation on graphs. 
\tm\label{thm:accdiscrete} There is a universal constant $C$ such that for any ancient solution $D_tu=\Delta u$ on $V\times \Z_-$ and $R\in \N, R\geq s,$
\begin{equation}\label{eq:accdiscrete1} R^2\sum_{\wt{Q}_R}\Gamma(u)+R^4\sum_{\wt{Q}_R}(D_tu)^2\leq C\sum_{\wt{Q}_{9R}} u^2.
\end{equation}
\tmd
\begin{proof} For any $R>0,$ let $\eta=\eta_R,$ where $\eta_R$ is the cut-off function defined in \eqref{eq:cut1}.

We first estimate $\sum_{\wt{Q}_R}\Gamma(u).$ For any $t\in \Z_-,$ taking the time difference of $\sum\eta^2u^2$ at $t,$ we have, by $D_tu=\Delta u$
\begin{equation*}
D_t(\sum\eta^2 u^2)=\sum\eta^2 D_t(u^2)\leq \sum 2\eta^2 uD_tu=2\sum\eta^2u\Delta u,
\end{equation*} where we have used Proposition~\ref{prop:ineq1}. This is a similar result for \eqref{eq:energy1} in the proof of Theorem~\ref{thm:paracacc}. Applying the same argument therein, we get, for any $R\geq s,$
\begin{equation*}
D_t(\sum\eta^2 u^2)\leq-\sum_{}\Gamma(u)\eta^2+\frac{2}{R^2}\sum_{B_{3R}}u^2,
\end{equation*} which is an analog of \eqref{eq:bas11}.

Fix $R\in \N, R\geq s.$ For any $T\in \N,$ by summing over $t$ from $-T$ to $0$ in the above inequality, we obtain,  by Proposition~\ref{prop:elem1}, 
\begin{eqnarray}
\sum_{t=-T}^0\sum_{B_R}\Gamma(u)&\leq& \frac{2}{R^2}\sum_{t=-T}^0\sum_{B_{3R}}u^2+\left.\sum\eta^2 u^2\right|_{t=-T-1}\nonumber\\
&\leq &\frac{2}{R^2}\sum_{t=-T}^0\sum_{B_{3R}}u^2+\left.\sum_{B_{2R}} u^2\right|_{t=-T-1}.\label{eq:estdis1}
\end{eqnarray}
By Proposition~\ref{prop:elem2}, there exists $T_3\in [R^2+1,4R^2]\cap \N$ such that
$$\left.\sum_{B_{2R}}u^2\right|_{t=-T_3}=\frac{1}{3R^2}\sum_{t=-4R^2}^{-R^2-1}\sum_{x\in B_{2R}}m_xu^2(x,t).$$ By applying \eqref{eq:estdis1} for $T=T_3-1$ and the above equation, we get, for $R\in \N, R\geq s,$
\begin{eqnarray}
\sum_{\wt{Q}_R}\Gamma(u)&\leq& \sum_{t=-T_3+1}^0\sum_{B_R}\Gamma(u)
\leq \frac{2}{R^2}\sum_{t=-T_3+1}^0\sum_{B_{3R}}u^2+\left.\sum_{B_{2R}} u^2\right|_{t=-T_3}\nonumber\\
&\leq &\frac{C}{R^2}\sum_{\wt{Q}_{3R}}u^2+\frac{1}{3R^2}\sum_{t=-4R^2}^{-R^2-1}\sum_{B_{2R}}u^2\nonumber\\
&\leq&\frac{C}{R^2}\sum_{\wt{Q}_{3R}}u^2.\label{eq:estdis2}
\end{eqnarray}

Now we estimate $\sum_{\wt{Q}_R}(D_tu)^2.$ Let $h(t)$ be the function defined in \eqref{eq:defht}.  For any $t\in \Z_-,$ taking the time difference of $h(t)$ at $t,$ we get
\begin{equation*}D_t h\leq\sum_{x,y}w_{xy}\nabla_{xy}u\nabla_{xy}(D_t u)\eta(x)\eta(y),\end{equation*}
where we have used Proposition~\ref{prop:ineq1}. This is a similar result for \eqref{eq:energy2} in the proof of Theorem~\ref{thm:paracacc}. Applying the same argument therein, we get, for any $R\in \N, R\geq s,$ 
\begin{eqnarray*}D_th\leq-\sum (D_t u)^2\eta^2+\frac{2}{R^2}\sum_{B_{3R}}\Gamma(u).
\end{eqnarray*}

Fix $R\in \N, R\geq s.$ For any $T\in \N,$ by summing over $t$ from $-T$ to $0$ in the above inequality, we obtain,  by Proposition~\ref{prop:elem1}, 

\begin{eqnarray}
\sum_{t=-T}^0\sum_{B_R}(D_t u)^2&\leq&\sum_{t=-T}^0\sum (D_t u)^2\eta^2\leq \frac{2}{R^2}\sum_{t=-T}^0 \sum_{B_{3R}}\Gamma(u)+
h(-T-1)\nonumber\\
&\leq &\frac{2}{R^2}\sum_{t=-T}^0 \sum_{B_{3R}}\Gamma(u)+\left.\frac12\sum_{x,y\in B_{2R}}w_{xy}|\nabla_{xy}u|^2\right|_{t=-T-1}\nonumber\\
&\leq&\frac{2}{R^2}\sum_{t=-T}^0 \sum_{B_{3R}}\Gamma(u)+\left.\sum_{B_{2R}}\Gamma(u)\right|_{t=-T-1}.
\label{eq:etdis1}
\end{eqnarray}

By Proposition~\ref{prop:elem2}, there exists $T_4\in [R^2+1,4R^2]\cap \N$ such that
\begin{equation*}\left.\sum_{B_{2R}}\Gamma(u)\right|_{t=-T_4}=\frac{1}{3R^2}\sum_{t=-4R^2}^{-R^2-1}\sum_{x\in B_{2R}}m_x\Gamma(u)(x,t)\leq \frac{1}{3R^2}\sum_{\wt{Q}_{2R}}\Gamma(u).\end{equation*}
By applying \eqref{eq:etdis1} for $T=T_4-1,$ and by the above inequality, we get, for $R\in \N, R\geq s,$
\begin{eqnarray}
\sum_{\wt{Q}_R}(D_t u)^2&\leq&\sum_{t=-T_4+1}^0\sum_{B_R}(D_t u)^2\leq\frac{2}{R^2}\sum_{t=-T_4+1}^0 \sum_{B_{3R}}\Gamma(u)+\left.\sum_{B_{2R}}\Gamma(u)\right|_{t=-T_4}\nonumber\\
&\leq &\frac{C}{R^2}\sum_{\wt{Q}_{3R}}\Gamma(u)+\frac{1}{3R^2}\sum_{\wt{Q}_{2R}}\Gamma(u)\nonumber\\
&\leq&\frac{C}{R^2}\sum_{\wt{Q}_{3R}}\Gamma(u)\leq \frac{C}{R^4}\sum_{\wt{Q}_{9R}}u^2,
\end{eqnarray} where we have used \eqref{eq:estdis2} in the last inequality.

This proves the theorem.

\end{proof}

We recall some facts on difference operators on $\Z_+:=\Z\cap[0,+\infty).$
For any $f:\Z_+\to \R,$ we define $$\delta_n f(n_0)= f(n_0+1)-f(n_0),\ \forall n_0\in \Z_+.$$
Note that $\delta_n f(n_0)=-D_{t} g(-n_0),$ where $g(t)=f(-t).$
The binomial coefficients are defined as, for any $n\in \N,$ 
\[\binom{n}{i}=\left\{\begin{array}{ll}\frac{n(n-1)\cdots (n-i+1)}{i!},&0\leq i\leq n, i\in \Z,\\
0,& \mathrm{otherwise}.
\end{array}\right.\] Note that for any $i\in \N,$ $$\delta_n \binom{n}{i}=\binom{n+1}{i}-\binom{n}{i}=\binom{n}{i-1}.$$ Therefore, for the difference operators $\delta_n$ (or $D_t$),
$f(n)=\binom{n}{i}$ stands as an analog of the monomial $t^i,$ up to some factor, for the differential operator $\partial_t.$

The following proposition is well-known.
\begin{prop}\label{prop:bin} Let $f:\Z_+\to \R$ satisfy, for some $q\in \N,$ $$\delta_n^q f\equiv 0,\quad \mathrm{on}\ \Z_+,$$ where $\delta_n^q$ is the $q$-th composition of $\delta_n.$ 
Then there are $\{a_i\}_{i=1}^{q-1}\subset \R$ such that
$$f(n)=\sum_{i=0}^{q-1} a_i \binom{n}{i}.$$
\end{prop}

This yields the following corollary.
\co\label{coro:cdis1} Let $G$ be a weighted graph admitting an intrinsic metric satisfying Assumption~\ref{ass:1}. Suppose that $G$ has polynomial volume growth, i.e. \eqref{eq:polynomial1}, and $u\in \wt{\mathcal{P}}_{k}(G)$ for some $k>0.$  Then for any $q\in \N,$ $4q>2k+\alpha+2,$
$$D_t^q u\equiv 0.$$ In particular, there exist some functions $p_i(x),$ $1\leq i\leq q-1,$ such that
$$u(t,x)=\sum_{i=1}^{q-1}p_i(x)\binom{-t}{i}.$$
\cod
\begin{proof}
We follow the argument in \cite{ColdingM19}, see the proof of Corollary~\ref{coro:c1}. 

For the first assertion,
since $\Delta$ commutes with $D_t,$ for any $i\in\N,$ $D_t^i u$ is also an ancient solution of the discrete-time heat equation.  For any $R\in \N, R\geq s,$ by the same argument as in 
Corollary~\ref{coro:c1} and $4q>2k+\alpha+2,$ we get
\begin{eqnarray*}
\sum_{\wt{Q}_R}|D_t^qu|^2\leq \frac{C(q)}{R^{4q}}\sum_{\wt{Q}_{9^qR}} u^2
\leq C R^{-4q+2k+\alpha+2}\to 0,\quad R\to\infty.
\end{eqnarray*} This proves the first assertion.

The second assertion follows from the first one and Proposition~\ref{prop:bin}.
\end{proof}

Now we can prove Theorem~\ref{thm:main2}. 
\begin{proof}[Proof of Theorem~\ref{thm:main2}]
The proof follows verbatim from \cite{ColdingM19}. 
Choose $q\in \N$ such that $4q>4k+\alpha+2.$ By Corollary~\ref{coro:cdis1}, we have
$$u(x,t)=p_0(x)+p_1(x)\binom{-t}{1}+\cdots+p_{q-1}(x)\binom{-t}{q-1},\quad x\in V, t\in \Z_-.$$ 

Note that $u\in \wt{\mathcal{P}}_{2k}(G).$ For any fixed $x\in V,$ considering sufficiently negative $t\in \Z_-$ in the above equality, we obtain that $$p_i(x)=0, \quad \forall i>k.$$ This yields that
\begin{equation}\label{eq:dpp1}u(x,t)=p_0(x)+p_1(x)\binom{-t}{1}+\cdots+p_{l}(x)\binom{-t}{l},
\end{equation} where $l:=\lfloor k\rfloor,$ the greatest integer less than or equal to $k.$

We claim that the function $p_i(x),$ $0\leq i\leq l,$ grows polynomially with the growth rate less than or equal to $2k.$ Fix distinct integer values $t_{l+1}<\cdots<t_2<t_1<-l.$ Set column vectors
$$\beta_j:=\left(1,\binom{-t_j}{1},\binom{-t_j}{2},\cdots, \binom{-t_j}{l}\right)^T, \quad 1\leq j\leq l+1.$$ Let 
$$B:=(\beta_1,\beta_2,\cdots,\beta_{l+1}).$$ Note that \[\det B =\frac{1}{\prod_{i=1}^l(i!)}
\det\begin{pmatrix}
    1 & -t_1 & (-t_1)^2 & \dots  & (-t_1)^l \\
   1 & -t_2 & (-t_2)^2 & \dots  & (-t_2)^l \\
    \vdots & \vdots & \vdots & \ddots & \vdots \\
    1 & -t_{l+1} & (-t_{l+1})^2 & \dots  & (-t_{l+1})^l \\
\end{pmatrix}\neq 0,
\] which yields that $\{\beta_j\}_{j=1}^{l+1}$ are linear independent in $\R^{l+1}.$ Hence there are $b_{j}^i\in \R$ such that
$$p_i(x)=\sum_{j=1}^{l+1}b_j^i u(x,t_j).$$ This proves the claim.

Since $D_t u=\Delta u,$ by \eqref{eq:dpp1}, 
\begin{equation*} \Delta p_l=0,\ \Delta p_i=-p_{i+1}, \quad 0\leq i\leq l-1.
\end{equation*} Then by the linear algebra argument as in the proof of Theorem~\ref{thm:main1}, we get
\begin{eqnarray*}\dim \wt{\mathcal{P}}_{2k}(G)&\leq &(k+1) \mathcal{H}_{2k}(G).
\end{eqnarray*} This proves the theorem.

\end{proof}

\textbf{Acknowledgements.} We thank Qi S. Zhang for many discussions and comments on ancient solutions of the heat equation. The author is supported by NSFC, no.11831004 and no. 11926313.


\bibliographystyle{alpha}
\bibliography{cal-poly-unbounded}

\end{document}